
\def\LaTeX{%
  \let\Begin\begin
  \let\End\end
  \let\salta\relax
  \let\finqui\relax
  \let\futuro\relax}

\def\UK{\def\our{our}\let\sz s}
\def\USA{\def\our{or}\let\sz z}



\LaTeX

\USA


\salta

\documentclass[twoside,12pt]{article}
\setlength{\textheight}{24cm}
\setlength{\textwidth}{16cm}
\setlength{\oddsidemargin}{2mm}
\setlength{\evensidemargin}{2mm}
\setlength{\topmargin}{-15mm}
\parskip2mm


\usepackage{amsmath}
\usepackage{amsthm}
\usepackage{amssymb}
\usepackage[mathcal]{euscript}

\usepackage[usenames,dvipsnames]{color}
%
%


\def\gianni{\color{green}}
\let\gianni\relax
\def\pier{\color{magenta}}
\let\pier\relax
\def\juerg{\color{red}}
\def\pcol{\color{blue}}
\let\juerg\relax
\let\pcol\relax


\bibliographystyle{plain}


%

\finqui

\def\Beq{\Begin{equation}}
\def\Eeq{\End{equation}}
\def\Bsist{\Begin{eqnarray}}
\def\Esist{\End{eqnarray}}

\def\Bthm{\Begin{theorem}}
\def\Ethm{\End{theorem}}

\def\Bprop{\Begin{proposition}}
\def\Eprop{\End{proposition}}

\def\Brem{\Begin{remark}\rm}
\def\Erem{\End{remark}}

\def\Bdim{\Begin{proof}}
\def\Edim{\End{proof}}
\let\non\nonumber




\def\step #1 \par{\medskip\noindent{\bf #1.}\quad}


\def\holder{H\"older}
\def\aand{\quad\hbox{and}\quad}

\def\rhs{right-hand side}


\def\generaliz{generali\sz}

\def\organiz{organi\sz}


\def\multibold #1{\def\arg{#1}%
  \ifx\arg\pto \let\next\relax
  \else
  \def\next{\expandafter
    \def\csname #1#1#1\endcsname{{\bf #1}}%
    \multibold}%
  \fi \next}

\def\pto{.}

\def\multical #1{\def\arg{#1}%
  \ifx\arg\pto \let\next\relax
  \else
  \def\next{\expandafter
    \def\csname cal#1\endcsname{{\cal #1}}%
    \multical}%
  \fi \next}


\def\multimathop #1 {\def\arg{#1}%
  \ifx\arg\pto \let\next\relax
  \else
  \def\next{\expandafter
    \def\csname #1\endcsname{\mathop{\rm #1}\nolimits}%
    \multimathop}%
  \fi \next}

\multibold
qwertyuiopasdfghjklzxcvbnmQWERTYUIOPASDFGHJKLZXCVBNM.

\multical
QWERTYUIOPASDFGHJKLZXCVBNM.

\multimathop
dist div dom meas sign supp .


\def\accorpa #1#2{\eqref{#1}--\eqref{#2}}
\def\Accorpa #1#2 #3 {\gdef #1{\eqref{#2}--\eqref{#3}}%
  \wlog{}\wlog{\string #1 -> #2 - #3}\wlog{}}


\def\<#1>{\mathopen\langle #1\mathclose\rangle}
\def\norma #1{\mathopen \| #1\mathclose \|}

\def\iot {\int_0^t}

\def\iO{\int_\Omega}
\def\intQt{\int_{Q_t}}
\def\intQs{\int_{Q_s}}

\def\smiot{\mathop{\textstyle\iot}\nolimits}

\def\dt{\partial_t}
\def\dn{\partial_\nu}
\def\ds{\,ds}

\def\cpto{\,\cdot\,}

\def\checkmmode #1{\relax\ifmmode\hbox{#1}\else{#1}\fi}
\def\aeO{\checkmmode{a.e.\ in~$\Omega$}}
\def\aeQ{\checkmmode{a.e.\ in~$Q$}}


\def\erre{{\mathbb{R}}}




\def\genspazio #1#2#3#4#5{#1^{#2}(#5,#4;#3)}
\def\spazio #1#2#3{\genspazio {#1}{#2}{#3}T0}

\def\L {\spazio L}

\def\W {\spazio W}

\def\Vp{V^*}


\def\Lx #1{L^{#1}(\Omega)}
\def\Hx #1{H^{#1}(\Omega)}
\def\Wx #1{W^{#1}(\Omega)}

\def\Ldue{\Lx 2}

\def\Huno{\Hx 1}


\def\LQ #1{L^{#1}(Q)}


\let\theta\vartheta
\let\eps\varepsilon

\let\TeXchi\chi                         
\newbox\chibox
\setbox0 \hbox{\mathsurround0pt $\TeXchi$}
\setbox\chibox \hbox{\raise\dp0 \box 0 }
\def\chi{\copy\chibox}


\def\muz{\mu_0}
\def\rhoz{\rho_0}

\def\rhomin{\rho_*}
\def\rhomax{\rho^*}
\def\ximin{\xi_*}
\def\ximax{\xi^*}

\def\kamu{\kappa(\mu)}
\def\kamut{\kappa(\mu(t))}

\def\Kmu{K(\mu)}

\def\Kmus{K(\mu(s))}
\def\ktilde{\widetilde k}
\def\kmin{\kappa_*}
\def\kmax{\kappa^*}
\def\ximin{\xi_*}
\def\ximax{\xi^*}




\def\Krejci{Krej\v c\'\i{}}

\Begin{document}

%

%


\title{\bf A {\pier continuous dependence} result
for a nonstandard system of phase field equations}
\author{}
\date{}
\maketitle
\begin{center}
\vskip-2.6cm
{\large\bf Pierluigi Colli$^{(1)}$}\\
{\normalsize e-mail: {\tt pierluigi.colli@unipv.it}}\\[.25cm]
{\large\bf Gianni Gilardi$^{(1)}$}\\
{\normalsize e-mail: {\tt gianni.gilardi@unipv.it}}\\[.25cm]
{\large\bf Pavel \Krejci $^{(2)}$}\\
{\normalsize e-mail: {\tt krejci@math.cas.cz}}\\[.25cm]
{\large\bf J\"urgen Sprekels$^{(3)}$}\\
{\normalsize e-mail: {\tt sprekels@wias-berlin.de}}\\[.45cm]
$^{(1)}$
{\small Dipartimento di Matematica ``F. Casorati'', Universit\`a di Pavia}\\
{\small via Ferrata 1, 27100 Pavia, Italy}\\[.2cm]
$^{(2)}$
{\small Institute of Mathematics, Czech Academy of Sciences}\\
{\small \v{Z}itn\'a 25, CZ-11567 Praha 1, Czech Republic}\\[.2cm]
$^{(3)}$
{\small Weierstra\ss-Institut f\"ur Angewandte Analysis und Stochastik}\\
{\small Mohrenstra\ss e\ 39, 10117 Berlin, Germany}\\[.5cm]
\end{center}
\vskip.1cm


\Begin{abstract}
\vskip-.5cm
Th{\pcol e present} note deals with a nonstandard systems of differential equations
{\pier describing a} two-species phase segregation.
{\juerg This system} naturally arises in the asymptotic analysis carried out 
recently by {\pier the same authors,} 
as the diffusion coefficient in the equation governing 
the evolution of the order parameter tends to zero.
In particular, an existence result {\pier has been} proved 
for the limit system in a very general framework.
On the contrary, uniqueness {\pier was shown} by 
assuming {\pier a {\sl constant}} mobility coefficient.
{Here, we \generaliz e} {\juerg this} result and {\juerg prove} 
a continuous dependence property in the case {\pier that the mobility 
coefficient suitably depends} on the chemical potential.
\\[0.3cm]
\End{abstract}

\bigskip


\salta

\pagestyle{myheadings}
\newcommand\testopari{\sc Colli \ --- \ Gilardi \ --- \ \Krejci{} \ --- \ Sprekels}
\newcommand\testodispari{\sc Continuous dependence for a nonstandard phase field system}
\markboth{\testodispari}{\testopari}

\finqui


\section{Introduction}
\label{Intro}
\setcounter{equation}{0}
In this paper, we {\pier address the system} 
\Bsist
  && \bigl( 1 + 2g(\rho) \bigr) \, \dt\mu
  + \mu \, g'(\rho) \, \dt\rho
  - \div \bigl( \kamu\nabla\mu \bigr) = 0
  \label{Iprima}
  \\
  && \dt\rho + f'(\rho) = \mu \, g'(\rho)
  \label{Iseconda}
  \\
  && (\kamu\nabla\mu) \cdot \nu|_\Gamma = 0
  \label{Ibc}
  \\
  && \mu(0) = \muz
  \aand
  \rho(0) = \rhoz 
  \label{Icauchy}
\Esist
\Accorpa\Ipbl Iprima Icauchy
{\pier of differential equations and boundary and initial conditions
in terms of} the unknown fields $\mu$ and~$\rho$; equations \accorpa{Iprima}{Iseconda} are meant to hold in a bounded domain $\Omega\subset\erre^3$ 
with a smooth boundary~$\Gamma$ and in some time interval~$(0,T)$,
and $\nu$ in \eqref{Ibc} denotes the {\juerg outward} unit normal vector {\juerg to} $\Gamma$. 
{\pier The recent paper~\cite{CGKSasy} investigate{\juerg d}
the existence of solutions to the above system: actually, a solution {\juerg was}
found by considering the analogous system in which} the ordinary differential equation~\eqref{Iseconda} is replaced by the partial differential equation
\Beq
  \dt\rho - \sigma \Delta\rho + f'(\rho) = \mu \, g'(\rho)
  \quad \hbox{with the boundary condition} \quad
  \dn\rho|_\Gamma = 0,
  \label{Isecondas}
\Eeq
and {\pier then} performing the asymptotic analysis as {\juerg the} {\pier diffusive 
coefficient} $\sigma$ tends to zero.

Such a modified system {\pier arises} from the model introduced in~\cite{Podio},
{\juerg which describes} the phase segregation of two species (atoms and vacancies, say) 
on~a lattice in {\juerg the} presence of diffusion. {\pier It turns out to be} a modification 
of the well-known Cahn-Hilliard equations (see, e.g.,~{\pcol\cite{FG, Gurtin}}).
The~state variables are the  {\sl order parameter\/}~$\rho$
(volume density of one of the two species),
which must of course {\juerg attain} values in the 
{\pier domain of {\juerg the} nonlinearities $g' $ and $f'$},
and the {\sl chemical potential\/}~$\mu$,
which is required to be nonnegative for physical reasons. 
{\pier The initial-boundary value problem for the PDE system} 
has been studied in {\pier a} series of papers {\pier with 
a number of obtained results: here, we confine ourselves to quote
the former \cite{CGPS3, CGPS7, CGPS6} and latter {\pcol\cite{CGKSasy,CGS2,CGKPS}}.} 

In the {\pier mentioned} papers, {\pier the function $g$ is taken as a smooth nonnegative and possibly concave function (like it looks here)}, while 
the function $f$ represents a {\pier multi-well} potential: in this respect, 
a thermodynamically relevant example {\juerg for} $f$ is the so-called {\sl logarithmic 
potential\/}{\pier , in which $f'$ is given} 
by~the formula
\Beq
  f'(\rho) = \ln \frac {1+\rho} {1-\rho} - 2c \rho 
  \quad \hbox{for $\rho\in(-1,1)$},
  \label{logpot} 
\Eeq
{\pier with} $c>1$ in order that $f$ actually presents a double well.
{\pier The class of the admissible potentials 
may be rather wide} and include both the 
standard double-well potential defined~by
\Beq
  f(\rho) = \frac 14 \, (\rho^2-1)^2
  \quad \hbox{for $\rho\in\erre$}
  \label{standardpot}
\Eeq
and potentials whose convex part $f_1$ is just a proper 
and lower semicontinuous function, thus possibly non-differentiable 
in its effective domain. In {\pier such a case},
the~monotone part $f_1'$ of $f'$ is replaced by 
the (possibly) multivalued subdifferential $\partial f_1$
and \eqref{Isecondas} has to be read as a differential inclusion. 
In~\cite{CGKSasy}, {\juerg this wide class of potentials was} considered.
Moreover, {\pier in \cite{CGKSasy} the mobility coefficient $\kappa$ in \eqref{Iprima} and \eqref{Ibc} 
{\juerg was} allowed to depend also on~$\rho$}.

{\pier Therefore, the existence result for system \Ipbl\ {\juerg proved} in 
\cite{CGKSasy}} is very general. {\pier On the other hand}, the solution constructed 
in this way is rather irregular, in principle (due~to a lack of regularity 
for~$\mu$).
Nevertheless, {\pier it has been shown to be unique (and~a little smoother than 
expected) provided that} the mobility coefficient $\kappa$ is a positive constant.
 
The {\pier aim of the present paper is {\juerg to generalize} the uniqueness proof 
performed in~\cite{CGKSasy} to the case of a mobility coefficient depending 
on the chemical potential, exactly as in~\eqref{Iprima} and~\eqref{Ibc}. 
{\pier  Moreover, the continuous dependence of the solution on the initial 
data is shown in terms of suitable norms. Of course, in order to accomplish our 
program, a natural uniform parabolicity condition is required for $\kappa$.}

The paper is \organiz ed as follows. In the next section, we list our assumptions 
and {\pier rewrite problem \Ipbl\ in a precise form. In Section~\ref{Results}, 
we state} and prove our {\pier uniqueness and continuous dependence} result.


\section{Assumptions and notations}
\label{MainResults}
\setcounter{equation}{0}

We first introduce precise assumptions on the data 
for the mathematical problem under investigation. 
We assume $\Omega$ to be a bounded connected 
open set in $\erre^3$ with smooth boundary~$\Gamma$
and let  $T\in(0,+\infty)$ stand for a final time. 
We~set for brevity 
\Beq
  V := \Huno, \quad
  H := \Ldue , 
  \aand
  Q := \Omega \times (0,T) .
  \label{defspazi}
\Eeq
The symbol $\<\cpto,\cpto>$ denotes the duality product 
between~$\Vp$, the dual space of~$V$, and~$V$ itself.
For the nonlinearities we assume that there exist real constants 
$\kmin$, $\kmax$, $\rhomin$, $\rhomax$, $\ximin$, and $\ximax$
{\juerg such} that the combined conditions listed below hold.
\Bsist
  && \kappa:[0,+\infty) \to \erre
  \quad \hbox{is continuous}
  \label{hpkappa}
  \\
  && 0 < \kmin \leq \kappa(m) \leq \kmax
  \quad \hbox{for every $m\geq 0$}
  \label{hpparab}
  \\
  && f = f_1 + f_2 \,, \quad f_1:\erre \to [0,+\infty], \quad f_2 :\erre \to \erre
  \label{hpf}
  \\
  && \hbox{$f_1$ is convex, proper, l.s.c. and $f_2$ is a $C^2$ function}
  \label{hpfbis}
  \\
  && \beta := \partial f_1
  \aand
  \pi := f_2'
  \label{defbetapi}
  \\
  && \hbox{$g\in C^2(\erre)$, \quad $g(r)\geq0$ \  and \ $g''(r) \leq 0$ \ for  $r\in\erre$}
  \label{hpg}
  \\
  && \hbox{$\pi$, $g$, and $g'$ are Lipschitz continuous}
  \label{hplip}
  \\
  \noalign{\allowbreak}
  && \rhomin,\,\rhomax \in D(\beta) , \quad
  \ximin \in \beta (\rhomin) , 
  \aand
  \ximax \in \beta (\rhomax)
  \label{hpxirhomm}
  \\
  && \ximin + \pi (\rhomin) \leq 0 \leq \ximax + \pi (\rhomax) 
  \aand
  g'(\rhomin) \geq 0 \geq g'(\rhomax) .
  \label{hpsegnimm}
\Esist
\Accorpa\Hpstruttura hpkappa hpsegnimm
Notice that important potentials like \eqref{logpot} and \eqref{standardpot}
fit the above requirements with suitable choices of $g$ and of the constants.
For the initial data, we require~that
\Bsist
  && \muz \in V \cap \Lx\infty
  \aand
  \muz \geq 0
  \quad \aeO 
  \label{hpmuz}
  \\
  && \rhoz \in V 
  \aand
  \rhomin \leq \rhoz \leq \rhomax
  \quad \aeO .
  \label{hprhoz}
\Esist
\Accorpa\Hpdati hpmuz hprhoz
\Accorpa\Hptutto hpkappa hprhoz

Now, we recall the part that follows from the asymptotic analysis performed in~\cite{CGKSasy}
and {\juerg is} of interest for the present paper.

\Bthm
\label{Esistenza}
Assume that~\Hptutto\ hold. 
Then there exists at least {\juerg one} triplet $(\mu,\rho,\xi)$ {\juerg that satisfies}
\Bsist 
  && \mu \in \L\infty H \cap \L2V \cap \LQ\infty
  \aand
  \mu \geq 0 \quad \aeQ
  \label{regmu}
  \\
  && \rho \in \L\infty V , \quad
  \dt\rho \in \LQ\infty ,
  \aand
  \rhomin \leq \rho \leq \rhomax \quad \aeQ
  \qquad
  \label{regrho}
  \\
  && \xi \in \LQ\infty, \quad
  \xi \in \beta(\rho) 
  \aand
  \ximin \leq \xi \leq \ximax
  \quad \aeQ
  \label{regxi}
  \\
  && u := \bigl( 1+2g(\rho) \bigr) \mu \ \in \W{1,p}\Vp \cap \L2{\Wx{1,q}}
  \label{regu}
\Esist
\Accorpa\Regsoluz {regmu} {regu}
for some $p,\, q>1$ and {\juerg solves} the problem 
\Bsist
  && \< \dt u(t), v >
  + \iO  \kamut \nabla\mu(t) \cdot \nabla v
  = \iO \mu(t) \, g'(\rho(t)) \, \dt\rho(t) \, v 
  \non
  \\
  && \hskip5cm 
  \hbox{for all $v\in V$ and a.a.\ $t\in(0,T)$}
  \label{prima}
  \\
  && \dt\rho + \xi + \pi(\rho)
   = \mu \, g'(\rho)
  \quad \aeQ
  \label{seconda}
  \\
  && u(0) = \bigl( 1+2g(\rhoz) \bigr) \muz
  \aand
  \rho(0) = \rhoz
  \quad \aeO .
  \label{cauchy}
\Esist
\Accorpa\Pbl prima cauchy
\Ethm

\Brem
\label{Commentosuprima}
We notice that \eqref{prima} actually is a weak form of equation~\eqref{Iprima}
(with the boundary condition~\eqref{Ibc} since the test function $v$ is free on the boundary).
Indeed, whenever $\mu$ is smoother with respect to time,
one can compute $\dt u$ by the Leibniz rule
and see that the differential equation hidden in the variational equation \eqref{prima} 
coincides with~\eqref{Iprima}.
We also observe that \cite{CGKSasy} precisely yields $p=4/3$ and $q=3/2$ in~\eqref{regu}. {\pier On the other hand, the regularity 
$u\in \L2V $ follows immediately from  $u = \bigl( 1+2g(\rho) \bigr) \mu$ thanks to \accorpa{regmu}{regrho}, whence one can take $q=2$.}
\Erem

{\juerg In~\cite{CGKSasy}, an additional result was proved that} deals with continuous dependence on the initial datum~$\rhoz$.
{\pier Here, we adapt the statement to our purposes.}
{\gianni Indeed, we also consider possibly different initial data for the chemical potential
(even though they do not enter the final estimate, directly).}

\Bprop
\label{Perunicita}
{\pier Assume that~\Hptutto\ hold. Let 
{\gianni $(\mu_{0,i},\rho_{0,i})$, $i=1,2$, be two sets of initial data 
satisfying \eqref{hpmuz}--\eqref{hprhoz}},
and let $(\mu_i,\rho_i,\xi_i )$, $i=1,2$,  
be two solutions to the corresponding} problem \Pbl\ that
satisfy the regularity assumptions \Regsoluz.
Then the following estimate holds true{\juerg :} {\pier
\Bsist
  \hskip-1.5cm&& |(\rho_1 - \rho_2)(t)| + \iot 
   \Bigl(|\dt(\rho_1 -\rho_2)| + |\xi_1 - \xi_2|  \Bigr) (s)\ds    
  \non
  \\ 
  \hskip-1.5cm&& {}\leq{} C \left( | \rho _{0,1} -\rho_{0,2}| + \iot
   \bigl( |\mu_1 - \mu_2| + (1+\mu_1) |\rho_1 - \rho_2| \bigr)(s) \ds
 \right)  \label{perunicita}
\Esist
}
for every $t\in[0,T]$ and \aeO,
where $C$ depends only on the constants and the functions 
mentioned in our assumptions \Hpstruttura\ on the structure of the system.
\Eprop
\Bdim
{\pier In order to give just an idea {\juerg  how} to obtain \eqref{perunicita},
let us point out that the procedure consists in testing the difference of two equations \eqref{seconda} by $\sign (\xi_1 - \xi_2)$. Indeed, setting 
$w_i = \mu_i g'(\rho_i) - \pi(\rho_i)$, $i=1,2$, and multiplying the identity
\Beq
  \dt (\rho_1 - \rho_2) +  (\xi_1-\xi_2) = w_1 - w_2
  \label{prelip1}
\Eeq
by $\sign(\xi_1-\xi_2)$, it is not difficult to infer~that
\Beq
 \dt |\rho_1 - \rho_2| + | \xi_1-\xi_2|  \le |w_1 - w_2|
  \quad \hbox{a.e.\ in $Q$}.
  \label{lip1}
\Eeq
Thanks to the Lipschitz continuity  properties in \eqref{hplip},
and integrating \eqref{lip1} only with respect to time, we obtain {\juerg that for~$t\in(0,T)$ it holds}
\Bsist
   &&|\rho_1-\rho_2| (t) + \iot |\xi_1-\xi_2|(s) \, ds
   \non
   \\ 
    &&{}\leq c \left( | \rho _{0,1} -\rho_{0,2}| + \iot
   \bigl( |\mu_1 - \mu_2| + (1+\mu_1) |\rho_1 - \rho_2| \bigr)(s) \ds
 \right)
  \non
\Esist
a.e. in $\Omega$. Moreover, note that \eqref{prelip1} implies 
\Beq
  \iot \dt |\rho_1 - \rho_2|(s) \, ds 
  \leq \iot \bigl( |w_1-w_2| + |\xi_1 - \xi_2| \bigr)(s) \, ds,
  \non
\Eeq
whence \eqref{perunicita} {\juerg easily follows}.
}
\Edim

In \cite{CGKSasy}, {\juerg the uniqueness of the solution given by Theorem~\ref{Esistenza}
(as~well as the regularity $\dt\mu\in\LQ2$)
was proved under an additional assumption, namely:}

\Bthm
\label{Unicitapartic}
Assume \Hptutto\ and that $\kappa$ is a positive constant.
Then the solution $(\mu,\rho,\xi)$ given by Theorem~\ref{Esistenza} is unique.
\Ethm

The aim of this paper is to improve {\juerg this} result
by showing that uniqueness {\pier and continuous dependence hold} 
in the more general framework of Theorem~\ref{Esistenza},
as stated in the forthcoming Theorem~\futuro\ref{Unicita}.


\section{{\pier Uniqueness and continuous dependence}}
\label{Results}
\setcounter{equation}{0}

In this section, we prove the uniqueness {\pier and continuous dependence} result for the solution to problem \Pbl\ stated below.

\Bthm
\label{Unicita}
Assume that {\juerg the} conditions~\Hptutto\ are satisfied.
Then the solution $(\mu,\rho,\xi)$
given by Theorem~\ref{Esistenza} is unique. {\pier Moreover, let 
$(\mu_{0,i},\rho_{0,i})$, $i=1,2$, be two sets of initial data 
satisfying \eqref{hpmuz}--\eqref{hprhoz}, 
and let $(\mu_i,\rho_i,\xi_i )$, $i=1,2$,  
be the corresponding solutions, which fulfill \Pbl\ with 
{\gianni $\muz=\mu_{0,i}$ and $\rhoz=\rho_{0,i}$, $i=1,2$}.
Then there exists a constant $C$, depending on the data 
through the structural assumptions, such that}
\Bsist
\norma{\mu_1 -\mu_2}_{\L2 H} 
+\norma{\rho_1 -\rho_2}_{\L\infty H} + \norma{\xi_1 -\xi_2}_{L^1(Q)}
\non \\
\leq C \left\{ \norma{\mu_{0,1} -\mu_{0,2}}_H +  \norma{\rho_{0,1} -\rho_{0,2}}_H   \right\} . \label{pier}
\Esist
\Ethm
\Bdim
{\pier We just prove continuous dependence, with uniqueness as a byproduct}.
Throughout the proof,
we account for the well-known \holder\ inequality
and for the elementary Young inequality
\Beq
  ab\leq \eps a^2 + \frac 1 {4\eps} \, b^2
  \quad \hbox{for every $a,b\geq 0$ and $\eps>0$} .
  \non
\Eeq
Moreover, in order to simplify the notation, 
we use the same symbol small-case $c$ for different constants, 
{\juerg which} may only depend on~$\Omega$, the final time~$T$, the nonlinearities $\kappa$, $f$,~$g$,
and the solutions {\juerg under consideration. Thus,} 
the meaning of $c$ {\juerg may}
change from line to line and even {\juerg within} the same chain of inequalities.
{\gianni In contrast, we choose capital letters to denote precise constants we want to refer~to.}
Finally, we~set
\Beq
  Q_t := \Omega \times (0,t)
  \quad \hbox{for $t\in(0,T]$}.
  \label{defQt}
\Eeq
Our argument relies on {\pier a suitable adaptation of} the technique 
developed in~\cite{CGKSasy}
with the help of the function $K:[0,+\infty)\to\erre$ 
defined~by
\Beq
  K(m) := \int_0^m \kappa(m') \, dm'
  \quad \hbox{for $m\geq0$} .
  \label{defK}
\Eeq
We have indeed $K'=\kappa$, whence $\nabla\Kmu=\kamu\nabla\mu$,
so that \eqref{prima} becomes{\juerg , with $k:=\Kmu$,}
\Bsist
  \< \dt u(t) , v >
  + \iO \nabla k(t) \cdot \nabla v 
  = \iO \mu(t) \, g'(\rho(t)) \, \dt\rho(t) \, v
  \non
  \\
  \hskip5cm 
  \hbox{{\pcol for all $v\in V$ and a.a.\ $t\in(0,T)$.}}
  \label{primabis}
\Esist
We remark at once that \eqref{hpparab} yields {\juerg that for every $m_1\,,m_2\geq0$
it holds}
\Bsist
  && \kmin |m_1-m_2|
  \leq |K(m_1)-K(m_2)|
  \leq \kmax |m_1-m_2|
  \label{propKa}
  \\
  && (m_1-m_2) \bigl( K(m_1)-K(m_2) \bigr)
  \geq \kmin (m_1-m_2)^2 .
  \label{propKb}
\Esist
In our proof, we use the equation obtained {\juerg by} integrating \eqref{primabis} with respect to time
rather than \eqref{primabis} itself; namely, for every $v\in V$ and $t\in[0,T]$ {\juerg we have}
\Beq
  \iO u(t) \, v - \iO u(0) \, v
  + \iO \nabla\ktilde(t) \cdot \nabla v
  = \iO \Bigl( \smiot \mu(s) \, g'(\rho(s)) \, \dt\rho(s) \ds \Bigr) \, v,
  \label{intprima}
\Eeq
{\pier where}
\Beq
  \ktilde(t) := \iot k(s) \ds = \iot \Kmus \ds,
  \label{defktilde}
\Eeq
{\pier and {\juerg where}} $u(0)=(1+2g(\rhoz))\muz$ according to the first Cauchy condition~\eqref{cauchy}.
It is worth observing that the pointwise values of $u$ in the integrals over $\Omega$ are {\juerg well defined}.
Indeed, the boundedness of~$u$ (derived from the boundedness of $\mu$ and~$g(\rho)$) 
and the regularity of~$\dt u$
(cf.~\accorpa{regmu}{regrho} and~\eqref{regu})
ensure that $u$ is weakly continuous, e.g., as an $H$-valued function.
Now, we {\pier let 
$ (\mu_{0,i},  \rho_{0,i}) $, $i=1,2$, be the initial data and 
pick two solutions $(\mu_i,\rho_i,\xi_i)$, $i=1,2$. Then, let us} define the corresponding functions $u_i$, $k_i$, and~$\ktilde_i$
(according to \eqref{regu}, \eqref{primabis}, and~\eqref{defktilde}),
as well as the new ones $\gamma_i$, as~follows:
\Bsist
  && \gamma_i := 1 + 2g(\rho_i) \,, \quad
  u_i := \gamma_i \, \mu_i \,, \quad
  k_i := K(\mu_i)\,,  
  \non
  \\
  && \ktilde_i(t) := \iot k_i(s) \ds
  \quad \hbox{for $t\in[0,T]$}, \quad i=1,2,
  \non
\Esist
in order to simplify the notation.
For the same reason, we~set
\Bsist
  && {\pier \muz := \mu_{0,1} - \mu_{0,2} \,, \quad
  \rhoz := \rho_{0,1} - \rho_{0,2}}
  \non
  \\
  && \mu := \mu_1 - \mu_2 \,, \quad
  \rho := \rho_1 - \rho_2 \,, \quad
  \xi := \xi_1 - \xi_2 
  \non
  \\
  && \gamma := \gamma_1 - \gamma_2 \,, \quad
  u := u_1 - u_2 \,, \quad
  k := k_1 - k_2 \,, \aand
  \ktilde := \ktilde_1 - \ktilde_2 \,.
  \non
\Esist
At this point, we write \eqref{intprima} at the time $t$ for both solutions
and choose $v=k(t)=\dt\ktilde(t)$ in the difference.
We obtain
\Bsist
  && \iO u(t) k(t)
  + \iO \nabla\ktilde(t) \cdot \nabla\dt\ktilde(t) 
  \non
  \\
  && {\pier = \iO \bigl( (1+2 g(  \rho_{0,1}))\muz + 2 \mu_{0,2} ( g(  \rho_{0,1} ) - g(  \rho_{0,2})) \bigr) k(t) }
  \non
  \\
  && \quad{}+ \iO \Bigl(
    \smiot \bigl( \mu_1(s) g'(\rho_1(s)) \dt\rho_1(s) - \mu_2(s) g'(\rho_2(s)) \dt\rho_2(s) \bigr) \ds
  \Bigr) k(t) .\label{testprima}
\Esist
We estimate each term of \eqref{testprima} separately.
By accounting for \accorpa{hpg}{hplip}, \accorpa{regmu}{regrho}, {\juerg as well as for} \accorpa{propKa}{propKb},
we have
\Bsist
  && uk 
  = (\gamma_1 \mu_1 - \gamma_2 \mu_2) k
  = \gamma_1 \mu k + \gamma \mu_2 k
 \geq \kmin |\mu|^2 - c |\rho| \, |\mu|
  \quad \aeQ , 
  \quad \hbox{whence} \quad
  \non
  \\
  \noalign{\smallskip}
  && \iO u(t) k(t)
  \geq \kmin \iO |\mu(t)|^2 - c \iO |\rho(t)| \, |\mu(t)| 
  \geq {\pier \frac {3 \kmin}4}\iO |\mu(t)|^2 - c \iO |\rho(t)|^2 .
  \non
\Esist
Next, we clearly see that
\Beq
  \iO \nabla\ktilde(t) \cdot \nabla\dt\ktilde(t) 
  = \frac 12 \, \frac d {dt} \iO |\nabla\ktilde(t)|^2 .
  \non
\Eeq
{\pier With the help of \eqref{hplip}, \eqref{hpmuz}--\eqref{hprhoz} and \eqref{propKa}, we 
{\juerg can} control the first term on the \rhs\ of~\eqref{testprima}: 
\Bsist
&&\Bigl| \iO \bigl( (1+2 g(  \rho_{0,1}))\muz + 2 \mu_{0,2} ( g(  \rho_{0,1} ) - g(  \rho_{0,2})) \bigr) k(t) \Bigr| 
   \non \\
&&\leq  c  \iO \bigl(| \muz| + | \rhoz | \bigr) |\mu(t) |
\leq c {}\bigl( \norma{\muz}_H^2  + 
 \norma{\rhoz}_H^2\bigr)  + \frac \kmin 4 \iO |\mu(t)|^2  .
\non
\Esist
In order to estimate the second term,} we observe that
\Bsist
  && |\mu_1 g'(\rho_1) \dt\rho_1 - \mu_2 g'(\rho_2) \dt\rho_2|
  \non
  \\
  && \leq |\mu| \, |g'(\rho_1)| \, |\dt\rho_1|
  + |\mu_2| \, |g'(\rho_1) - g'(\rho_2)| \, |\dt\rho_1|
  + |\mu_2| \, |g'(\rho_2)| \, |\dt\rho|
  \non
  \\
  && \leq c \bigl( |\mu| + |\rho| + |\dt\rho| \bigr)
  \quad \aeQ \,,
  \non
\Esist
thanks to our regularity assumptions on the solutions and on the structure
(cf.~\accorpa{regmu}{regrho} and~\eqref{hplip}).
By owing to Proposition~\ref{Perunicita}, we deduce that
\Bsist
  && \Bigl| \iot \bigl( \mu_1(s) g'(\rho_1(s)) \dt\rho_1(s) - \mu_2(s) g'(\rho_2(s)) \dt\rho_2(s) \bigr) \ds \Bigr|
  \non
  \\
  && \leq c \iot \bigl( |\mu(s)| + |\rho(s)| + |\dt\rho(s)| \bigr) \ds
  \leq c {}{\pier \Bigl( |\rhoz| + \smiot \bigl( |\mu(s)| + |\rho(s)| \bigr) \ds \Bigr)}
  \quad \aeO .
  \non
\Esist
Therefore, we have 
\Bsist
  && \iO \Bigl( \smiot \bigl( \mu_1(s) g'(\rho_1(s)) \dt\rho_1(s) - \mu_2(s) g'(\rho_2(s)) \dt\rho_2(s) \bigr) \ds \Bigr) k(t)
  \non
  \\
  && 
  \leq c \iO \Bigl( {\pier |\rhoz| +} \smiot ( |\mu(s)| + |\rho(s)| ) \ds \Bigr) \, |\mu(t)|
  \non
  \\
  && \leq \frac \kmin 4 \iO |\mu(t)|^2
  + c \iO \left\{{\pier |\rhoz|^2 +{}} \Bigl( \smiot |\mu(s)| \ds \Bigr)^2 + \Bigl( \smiot |\rho(s)| \ds \Bigr)^2 \right\}
  \non
  \\
  && \leq \frac \kmin 4 \iO |\mu(t)|^2
  + c \iO \left\{ {\pier |\rhoz|^2 +{}} \smiot |\mu(s)|^2 \ds + \smiot |\rho(s)|^2 \ds \right\}
  \non
  \\
  && = \frac \kmin 4 \iO |\mu(t)|^2 {\pier {}+ c\Vert\rhoz\Vert_H^2
  + c \intQt |\mu|^2 + c }\intQt |\rho|^2 .
  \non
\Esist
By combining the above equalities and inequalities with \eqref{testprima}, 
we {\pier infer that}
\Beq
  \frac \kmin 4 \iO |\mu(t)|^2
  + \frac 12 \, \frac d {dt} \iO |\nabla\ktilde(t)|^2 
  \leq {\pier c {}\bigl( \norma{\muz}_H^2  + 
 \norma{\rhoz}_H^2\bigr) +{}} c \iO |\rho(t)|^2 
  + c \intQt |\mu|^2 + {\pier c} \intQt |\rho|^2\,, 
  \non
\Eeq
and an integration with respect to time yields
\Bsist
  && \frac \kmin 4 \intQt |\mu|^2
  + \frac 12 \, \iO |\nabla\ktilde(t)|^2 
  \non \\
  &&
  \leq {\pier c {}\bigl( \norma{\muz}_H^2  + 
 \norma{\rhoz}_H^2\bigr) +{}} 
 c \intQt |\rho|^2
  + \iot \Bigl( \intQs |\mu|^2
  + \intQs |\rho|^2 \Bigr) \ds 
  \non
  \\
  && \leq {\pier c {}\bigl( \norma{\muz}_H^2  + 
 \norma{\rhoz}_H^2\bigr) +{}} 
c \iot \Bigl( \intQs |\mu|^2 \Bigr) \ds + c \intQt |\rho|^2 . 
 \label{stimaprima}
\Esist
{\pier Now, let us consider \eqref{perunicita}{\juerg . Squaring} and integrating over $\Omega$, then applying H\"older's inequality on the \rhs,
we easily obtain {\juerg that} 
\Beq  
\iO |\rho (t)|^2  \leq  D \norma{ \rhoz }_H^2 + 
D \intQt |\mu|^2 +  D \intQt |\rho|^2 
\label{pier1}
\Eeq  
for some positive constant $D$. {\juerg Moreover}, by integrating  \eqref{perunicita} over $\Omega$ and then squaring, we arrive at  
\Beq  
\left(  \intQt |\xi| \right)^2  \leq  D \norma{ \rhoz }_H^2 + 
D \intQt |\mu|^2 +  D \intQt |\rho|^2 \,,
\label{pier2}
\Eeq  
{\gianni where $D$ is the same constant as before, without loss of generality}.
Hence, we multiply \eqref{stimaprima} by 
$12D/\kmin$ and add it to \eqref{pier1} and \eqref{pier2}. This computation 
leads to 
\Bsist
  && D \norma{ \mu}_{L^2(Q_t) }^2 + \norma{ \rho (t) }_{H}^2 
+ \norma{ \xi}_{L^1(Q_t) }^2 
  \non \\
  &&
  \leq c \bigl( \norma{\muz}_H^2  +  \norma{\rhoz}_H^2\bigr) 
+  c \int_0^t \norma{ \mu }_{ L^2(Q_s)}^2 \ds 
+  c \int_0^t \norma{ \rho (s) }_{H}^2 \ds  .
 \label{stimasecon}
\Esist
At this point, it suffices to apply the Gronwall lemma to deduce a variation of 
\eqref{pier} with the squared norms. Therefore, \eqref{pier} is completely proved.
}\Edim

\Brem
Clearly, just a few of the assumptions \Hptutto\ are used in the above proof.
The whole set of hypotheses has been listed in the statement of Theorem~\ref{Unicita}
in order to ensure both the existence of a solution satisfying \Regsoluz\
and the validity of estimate~\eqref{perunicita}, 
according to Theorem~\ref{Esistenza} and Proposition~\ref{Perunicita}.
\Erem

{\pcol
\section*{Acknowledgements}

The authors gratefully acknowledge some financial support from the MIUR-PRIN Grant 2010A2TFX2 ``Calculus of variations'' for PC and~GG, the GA\v CR Grant P201/10/2315 
and RVO:~67985840 for~PK, the DFG Research Center {\sc Matheon} 
in Berlin for~JS.
} 


\vspace{3truemm}

\Begin{thebibliography}{10}

\bibitem{CGKPS} 
P. Colli, G. Gilardi, P. \Krejci, P. Podio-Guidugli, J. Sprekels,
Analysis of a time discretization scheme for
a nonstandard viscous Cahn-Hilliard system,
{\pcol preprint (2013), pp.~1-31.}

\bibitem{CGKSasy} 
P. Colli, G. Gilardi, P. \Krejci, J. Sprekels,
A vanishing diffusion limit 
in a nonstandard system of phase field equations,
preprint WIAS-Berlin n.~1758 (2012), pp.~1-19.

{\pcol
\bibitem{CGS2} 
P. Colli, G. Gilardi, J. Sprekels,
Regularity of the solution to a nonstandard
system of phase field equations,
preprint (2013), pp.~1-11.}

\bibitem{CGPS3} 
P. Colli, G. Gilardi, P. Podio-Guidugli, J. Sprekels,
Well-posedness and long-time behaviour for 
a nonstandard viscous Cahn-Hilliard system, 
{\it SIAM J. Appl. Math.} {\bf 71} (2011) 1849-1870.

\bibitem{CGPS7} 
P. Colli, G. Gilardi, P. Podio-Guidugli, J. Sprekels,
Global existence for a strongly coupled 
Cahn-Hilliard system with viscosity, 
Boll. Unione Mat. Ital. (9) {\bf 5} (2012) 495-513.

\bibitem{CGPS6} 
P. Colli, G. Gilardi, P. Podio-Guidugli, J. Sprekels,
Global existence and uniqueness for a singular/degenerate  
Cahn-Hilliard system with viscosity, preprint 
WIAS-Berlin n.~1713 (2012), pp.~1-28{\pcol , to appear 
in {\it J. Differential Equations.}}  

{\pcol 
\bibitem{FG} 
E. Fried and M.E. Gurtin, 
Continuum theory of thermally induced phase transitions based on an order 
parameter, {\it Phys. D} {\bf 68} (1993) 326-343.
}

\bibitem{Gurtin} 
M. Gurtin, Generalized Ginzburg-Landau and
Cahn-Hilliard equations based on a microforce balance,
{\it Phys.~D\/} {\bf 92} (1996) 178-192.

\bibitem{Podio}
P. Podio-Guidugli, 
Models of phase segregation and diffusion of atomic species on a lattice,
{\it Ric. Mat.} {\bf 55} (2006) 105-118.

\End{thebibliography}


\End{document}

\bye